\documentclass[submission, Proceedings]{SciPost}

\usepackage{enumerate}
\usepackage{hyperref}

\usepackage{amsmath,amsfonts}
\usepackage{amssymb}
\usepackage{amscd}

\hypersetup{
    colorlinks,
    linkcolor={red!50!black},
    citecolor={blue!50!black},
    urlcolor={blue!80!black}
}

\usepackage{doi}

\DeclareSymbolFont{usualmathcal}{OMS}{cmsy}{m}{n}
\DeclareSymbolFontAlphabet{\mathcal}{usualmathcal}

\begin{document}

%
\def\cosech{\rm cosech}
\def\sech{\rm sech}
\def\coth{\rm coth}
\def\tanh{\rm tanh}
\def\half{{1\over 2}}
\def\third{{1\over3}}
\def\fourth{{1\over4}}
\def\fifth{{1\over5}}
\def\sixth{{1\over6}}
\def\seventh{{1\over7}}
\def\eigth{{1\over8}}
\def\ninth{{1\over9}}
\def\tenth{{1\over10}}
\def\bN{\mathop{\bf N}}
\def\R{{\rm I\!R}}
\def\Eins{{\mathchoice {\rm 1\mskip-4mu l} {\rm 1\mskip-4mu l}
{\rm 1\mskip-4.5mu l} {\rm 1\mskip-5mu l}}}
\def\Z{{\mathchoice {\hbox{$\sf\textstyle Z\kern-0.4em Z$}}
{\hbox{$\sf\textstyle Z\kern-0.4em Z$}}
{\hbox{$\sf\scriptstyle Z\kern-0.3em Z$}}
{\hbox{$\sf\scriptscriptstyle Z\kern-0.2em Z$}}}}
\def\abs#1{\left| #1\right|}
\def\com#1#2{
        \left[#1, #2\right]}
\def\square{\kern1pt\vbox{\hrule height 1.2pt\hbox{\vrule width 1.2pt
   \hskip 3pt\vbox{\vskip 6pt}\hskip 3pt\vrule width 0.6pt}
   \hrule height 0.6pt}\kern1pt}
      \def\boxop{{\raise-.25ex\hbox{\square}}}
\def\contract{\makebox[1.2em][c]{
        \mbox{\rule{.6em}{.01truein}\rule{.01truein}{.6em}}}}
\def\ltap{\ \raisebox{-.4ex}{\rlap{$\sim$}} \raisebox{.4ex}{$<$}\ }
\def\gtap{\ \raisebox{-.4ex}{\rlap{$\sim$}} \raisebox{.4ex}{$>$}\ }
\def\mn{{\mu\nu}}
\def\rs{{\rho\sigma}}
\newcommand{\Det}{{\rm Det}}
\def\Tr{{\rm Tr}\,}
\def\tr{{\rm tr}\,}
\def\sumij{\sum_{i<j}}
\def\e{\,{\rm e}}
\def\non{\nonumber\\}
\def\br{{\bf r}}
\def\bp{{\bf p}}
\def\bx{{\bf x}}
\def\by{{\bf y}}
\def\brhat{{\bf \hat r}}
\def\bv{{\bf v}}
\def\ba{{\bf a}}
\def\bE{{\bf E}}
\def\bB{{\bf B}}
\def\bA{{\bf A}}
\def\pa{\partial}
\def\dA{\partial^2}
\def\ddx{{d\over dx}}
\def\ddt{{d\over dt}}
\def\der#1#2{{d #1\over d#2}}
\def\lie{\hbox{\it \$}} 
\def\partder#1#2{{\partial #1\over\partial #2}}
\def\secder#1#2#3{{\partial^2 #1\over\partial #2 \partial #3}}
%
\def\be{\begin{equation}}
\def\ee{\end{equation}\noindent}
\def\bear{\begin{eqnarray}}
\def\ear{\end{eqnarray}\noindent}
\def\bec{\blue\begin{equation}}
\def\eec{\end{equation}\black\noindent}
\def\bearc{\blue\begin{eqnarray}}
\def\earc{\end{eqnarray}\black\noindent}
\def\benn{\begin{enumerate}}
\def\enn{\end{enumerate}}
\def\veject{\vfill\eject}
\def\ven{\vfill\eject\noindent}
%
\def\eq#1{{eq. (\ref{#1})}}
\def\eqs#1#2{{eqs. (\ref{#1}) -- (\ref{#2})}}
%
\def\totint{\int_{-\infty}^{\infty}}
\def\posint{\int_0^{\infty}}
\def\negint{\int_{-\infty}^0}
\def\pint{{\dps\int}{dp_i\over {(2\pi)}^d}}
%
\newcommand{\GeV}{\mbox{GeV}}
\def\FFdual{F\cdot\tilde F}
\def\bra#1{\langle #1 |}
\def\ket#1{| #1 \rangle}
\def\braket#1#2{\langle {#1} \mid {#2} \rangle}
\def\vev#1{\langle #1 \rangle}
\def\rightvac{\mid 0\rangle}
\def\leftvac{\langle 0\mid}
\def\ihbar{{i\over\hbar}}
\def\ge{\hbox{$\gamma_1$}}
\def\gz{\hbox{$\gamma_2$}}
\def\gd{\hbox{$\gamma_3$}}
\def\go{\hbox{$\gamma_1$}}
\def\gt{\hbox{\$\gamma_2$}}
\def\gth{\hbox{$\gamma_3$}} 
\def\gf{\hbox{$\gamma_5\;$}}
\def\slash#1{#1\!\!\!\raise.15ex\hbox {/}}
\newcommand{\slD}{\,\raise.15ex\hbox{$/$}\kern-.27em\hbox{$\!\!\!D$}}
\newcommand{\slpartial}{\raise.15ex\hbox{$/$}\kern-.57em\hbox{$\partial$}}
\newcommand{\PP}{\cal P}
\newcommand{\G}{{\cal G}}
\newcommand{\nc}{\newcommand}
\newcommand{\Fkala}{F_{\kappa\lambda}}
\newcommand{\Fkanu}{F_{\kappa\nu}}
\newcommand{\Flaka}{F_{\lambda\kappa}}
\newcommand{\Flamu}{F_{\lambda\mu}}
\newcommand{\Fmunu}{F_{\mu\nu}}
\newcommand{\Fnumu}{F_{\nu\mu}}
\newcommand{\Fnuka}{F_{\nu\kappa}}
\newcommand{\Fmuka}{F_{\mu\kappa}}
\newcommand{\Fkalamu}{F_{\kappa\lambda\mu}}
\newcommand{\Flamunu}{F_{\lambda\mu\nu}}
\newcommand{\Flanumu}{F_{\lambda\nu\mu}}
\newcommand{\Fkamula}{F_{\kappa\mu\lambda}}
\newcommand{\Fkanumu}{F_{\kappa\nu\mu}}
\newcommand{\Fmulaka}{F_{\mu\lambda\kappa}}
\newcommand{\Fmulanu}{F_{\mu\lambda\nu}}
\newcommand{\Fmunuka}{F_{\mu\nu\kappa}}
\newcommand{\Fkalamunu}{F_{\kappa\lambda\mu\nu}}
\newcommand{\Flakanumu}{F_{\lambda\kappa\nu\mu}}

\newcommand{\pb}{\bar{p}}
\newcommand{\ph}{\hat{p}}
\newcommand{\gb}{\bar{g}}
\newcommand{\gh}{\hat{g}}
\newcommand{\zb}{\bar{z}}
\newcommand{\zh}{\hat{z}}
\newcommand{\wh}{\hat w}
\newcommand{\wb}{\bar w}
\newcommand{\p}{p\!\!\!/~}
\newcommand{\pbdash}{\bar{p} \!\!\!/~}
\newcommand{\q}{q\!\!\!/~}
\newcommand{\B}{\beta \!\!\!/~}
\newcommand{\tb}{\bar{t}}

\nc{\spa}[3]{\left\langle#1\,#3\right\rangle}
\nc{\spb}[3]{\left[#1\,#3\right]}
\nc{\ksl}{\not{\hbox{\kern-2.3pt $k$}}}
\nc{\hf}{\textstyle{1\over2}}
\nc{\pol}{\varepsilon}
\nc{\tq}{{\tilde q}}
\nc{\esl}{\not{\hbox{\kern-2.3pt $\pol$}}}
\newcommand{\cL}{\cal L}
\newcommand{\D}{\cal D}
\newcommand{\Dhalf}{{D\over 2}}
\def\eps{\epsilon}
\def\epshalf{{\epsilon\over 2}}
\def\lag{( -\partial^2 + V)}
\def\freeexp{{\rm e}^{-\int_0^Td\tau {1\over 4}\dot x^2}}
\def\kinb{{1\over 4}\dot x^2}
\def\kinf{{1\over 2}\psi\dot\psi}
\def\expk{{\rm exp}\biggl[\,\sum_{i<j=1}^4 G_{Bij}k_i\cdot k_j\biggr]}
\def\expp{{\rm exp}\biggl[\,\sum_{i<j=1}^4 G_{Bij}p_i\cdot p_j\biggr]}
\def\expshort{{\e}^{\half G_{Bij}k_i\cdot k_j}}
\def\expabb{{\e}^{(\cdot )}}
\def\epseps#1#2{\varepsilon_{#1}\cdot \varepsilon_{#2}}
\def\epsk#1#2{\varepsilon_{#1}\cdot k_{#2}}
\def\kk#1#2{k_{#1}\cdot k_{#2}}
\def\G#1#2{G_{B#1#2}}
\def\Gp#1#2{{\dot G_{B#1#2}}}
\def\GF#1#2{G_{F#1#2}}
\def\Dab{{(x_a-x_b)}}
\def\Dsq{{({(x_a-x_b)}^2)}}
\def\PITD{{(4\pi T)}^{-{D\over 2}}}
\def\4piTD{{(4\pi T)}^{-{D\over 2}}}
\def\4piT4{{(4\pi T)}^{-2}}
\def\TintmD{{\dps\int_{0}^{\infty}}{dT\over T}\,e^{-m^2T}
    {(4\pi T)}^{-{D\over 2}}}
\def\Tintm4{{\dps\int_{0}^{\infty}}{dT\over T}\,e^{-m^2T}
    {(4\pi T)}^{-2}}
\def\Tintm{{\dps\int_{0}^{\infty}}{dT\over T}\,e^{-m^2T}}
\def\Tint{{\dps\int_{0}^{\infty}}{dT\over T}}
\def\np{n_{+}}
\def\nm{n_{-}}
\def\Np{N_{+}}
\def\Nm{N_{-}}
\newcommand{\slG}{{{\dot G}\!\!\!\! \raise.15ex\hbox {/}}}
\newcommand{\Gd}{{\dot G}}
\newcommand{\Gund}{{\underline{\dot G}}}
\newcommand{\Gdd}{{\ddot G}}
\def\GBd12{{\dot G}_{B12}}
\def\Dx{\dps\int{\cal D}x}
\def\Dy{\dps\int{\cal D}y}
\def\Dpsi{\dps\int{\cal D}\psi}
\def\dint#1{\int\!\!\!\!\!\int\limits_{\!\!#1}}
\def\ddtau{{d\over d\tau}}
\def\ie{\hbox{$\textstyle{\int_1}$}}
\def\iz{\hbox{$\textstyle{\int_2}$}}
\def\id{\hbox{$\textstyle{\int_3}$}}
\def\ldop{\hbox{$\lbrace\mskip -4.5mu\mid$}}
\def\rdop{\hbox{$\mid\mskip -4.3mu\rbrace$}}
%
\newcommand{\1}{{\'\i}}
\newcommand{\no}{\noindent}
\def\non{\nonumber}
\def\dps{\displaystyle}
\def\sy{\scriptscriptstyle}
\def\sy{\scriptscriptstyle}

\begin{center}{\Large \textbf{
Group invariants for Feynman diagrams
\\
}}\end{center}

\begin{center}
Idrish Huet\textsuperscript{1},
Michel Rausch de Traubenberg\textsuperscript{2} and
Christian Schubert\textsuperscript{3$\star$}
\end{center}

\begin{center}
{\bf 1} 
Facultad de Ciencias en F\'isica y Matem\'aticas, Universidad Aut\'onoma de Chiapas\\
 Ciudad Universitaria, Tuxtla Guti\'errez 29050, Mexico 
\\
{\bf 2} 
Universit\'e de Strasbourg, CNRS, IPHC UMR7178, F-67037 Strasbourg Cedex, France
\\
{\bf 3} 
Centro Internacional de Ciencias A.C. Campus UNAM-UAEM Avenida Universidad 1001 
Cuernavaca Morelos Mexico C.P. 62100
%
* christianschubert137@gmail.com
\end{center}

\begin{center}
\today
\end{center}


\definecolor{palegray}{gray}{0.95}
\begin{center}
\colorbox{palegray}{
  \begin{tabular}{rr}
  \begin{minipage}{0.1\textwidth}
    \includegraphics[width=20mm]{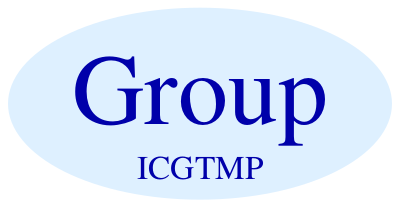}
  \end{minipage}
  &
  \begin{minipage}{0.85\textwidth}
    \begin{center}
    {\it 34th International Colloquium on Group Theoretical Methods in Physics}\\
    {\it Strasbourg, 18-22 July 2022} \\
    \doi{10.21468/SciPostPhysProc.?}\\
    \end{center}
  \end{minipage}
\end{tabular}
}
\end{center}

\section*{Abstract}
{\bf
It is well-known that the symmetry group of a Feynman diagram can give important information on possible strategies for its evaluation, and the mathematical objects that will be involved. 
Motivated by ongoing work on multi-loop multi-photon amplitudes in quantum electrodynamics, 
here I will discuss the usefulness of introducing a polynomial basis of invariants of the symmetry group of a diagram in Feynman-Schwinger parameter space. 
}

%

\section{Introduction: Schwinger parameter representation of Feynman diagrams}
\label{sec:intro}

The most universal approach to the calculation of Feynman diagrams uses Feynman-Schwinger parameters $ x_i$,
introduced through the exponentiation of the (Euclidean) scalar propagator,

$$ \frac{1}{p^2+m^2} = \int_0^\infty dx  \e^{-x (p^2+m^2)} $$

For scalar diagrams, one finds the following universal structure for an arbitrary graph $ G$ with $ n$ internal lines and $ l$ loops
in $ D$ dimensions:

$$ I_G = \Gamma(n-l D/2) \int_{x_i \geq 0} d^n x \, \delta\Bigl(1-\sum_{i=1}^n x_i\Bigr) \frac{{\cal U}^{n -  (l+1)\frac{D}{2}}}{{\cal F}^{n-l D/2}}$$

$ {\cal U}$ and $ {\cal F}$ are polynomials in the $ x_i$ called the first and second Symanzyk (graph) polynomials. There exist graphical methods for their construction. 

For more general theories (involving not only scalar particles) the same graph will involve, in addition to these two polynomials, also a 
{\it numerator polynomial} $ {\cal N}(x_1,\ldots,x_n)$. Such numerator polynomials can become extremely large, and in the present talk I would like to point
out the universal option of rewriting them as polynomials in a basis of invariants of the symmetry group of the graph.
After introducing the symmetries of Feynman diagrams and some related facts of invariant theory in section 2, sections 3 to 8
are devoted to the discussion and motivation of our main example, the three-loop effective Lagrangian in two-dimensional QED, 
where we have found this procedure to lead to significantly more manageable expressions. We come back to the group-theoretical
aspects of this calculation in section 9, before summarizing our findings in sections 10 and 11.

\section{Symmetries of Feynman diagrams}

Many Feynman diagrams possess a non-trivial symmetry group, generated by interchanges of the internal lines that leave the
topology of the graph unchanged 
\footnote{We do not consider here the exchange of external lines, since for such an exchange to leave a graph invariant
would require the two external momenta to be equal, which is not a natural condition. The exception is the
case where all the external momenta go to zero, which is what effectively happens in our main example below.}.
Then all its graph polynomials must be invariant under the natural action of the
group on the set of polynomials $\R[x_1,\ldots, x_n]$, $ g.P(X) \equiv P(g.X)$.

\begin{center}
\begin{figure}[htp]
\centerline{\includegraphics[width=5cm]{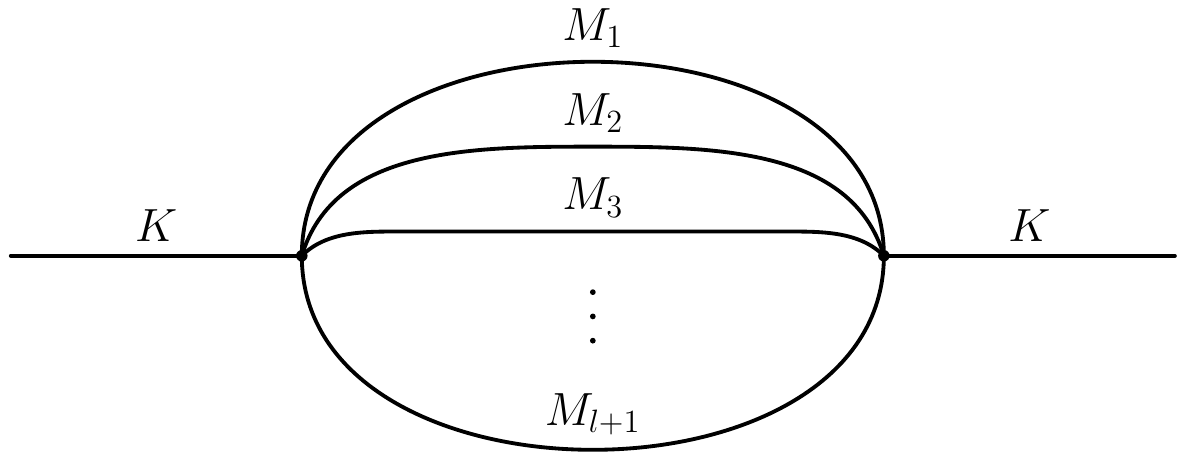}}
\vspace{-60pt}
\caption{$l-loop$ banana graph.}
\label{fig-banana}
\end{figure}
\end{center}

A nice example is the $ l-loop $ {\it banana graph} shown in Fig. \ref{fig-banana}.

\vspace{-80pt}

\vspace{-30pt}

When all the masses $M_i$ are equal,
this graph has full permutation symmetry in all the internal lines, so the symmetry group is $S_{l+1}$, and its graph polynomials must be symmetric functions of $x_1,\ldots, x_{l+1}$. As is well-known, this implies that they can be rewritten as
polynomials in the {\it elementary symmetric polynomials} $ S_1,\ldots,S_n$. 

Perhaps less known is that this generalizes to the case of a general symmetry group as follows \cite{sturm-book}:

{\bf Theorem 1}:
  Let $G$ be a finite group and let $\Gamma$ be an $n-$dimensional (real) representation.
  \begin{enumerate}
    \item There exist
  $n= \dim (\Gamma)$ algebraically invariant polynomials
  $P_1,\cdots, P_n$, called the {\it primitive invariants}, such that the Jacobian
$\frac {\partial (P_1,\ldots,P_n)}{\partial(x_1,\ldots,x_n)} \ne 0$.

\item
Denote $d_k=\deg(P_k)$ and ${\cal R} = \mathbb R[P_1,\cdots,P_n]$
the
subalgebra of  polynomial invariants generated by the primitive invariants.
\item There exist $m= d_1 \cdots d_n/|G|$
  secondary invariant polynomials $S_1,\cdots, S_m$.
  \item The subalgebra of invariants $\mathbb R[x_1,\cdots,x_n]^\Gamma$ is a free ${\cal R}-$module with basis $(S_1,\cdots,S_m)$.
In particular this means that any invariant $I \in  \mathbb R[x_1,\cdots,x_n]^\Gamma$ can be uniquely written as
$I=\sum\limits_{i=1}^m  f_i(P_1,\cdots,P_n) S_i$, 
where $f_i(P_1,\cdots,P_n), i=1,\cdots, m$ belong to ${\cal R}$, {\it i.e.}, are polynomials in $(P_1,\cdots,P_n)$.
\end{enumerate}

There exist computer algebra systems for the computation of $P_1,\ldots,P_n$ such as SINGULAR \cite{decker}.

\section{The Euler-Heisenberg Lagrangian at one loop}

\begin{figure}[ht]
\centerline{\includegraphics[scale=.6]{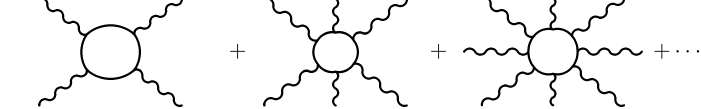}}
\vspace{10pt}
\caption{Feynman diagrams corresponding to the EHL.}
\label{fig-EHL1loop}
\end{figure}

In 1936 Heisenberg and Euler obtained their famous representation of the 
 one-loop QED effective Lagrangian in a constant field
(``Euler-Heisenberg Lagrangian'' or ``EHL'')
\bear
{\cal L}^{(1)}(a,b)&=& - {1\over 8\pi^2}
\int_0^{\infty}{dT\over T^3}
\,\e^{-m^2T}
\biggl\lbrack
{(eaT)(ebT)\over {\rm tanh}(eaT)\tan(ebT)} 
- {e^2\over 3}(a^2-b^2)T^2 -1
\biggr\rbrack
\label{ehspin}
\ear
Here $m$ is the electron mass, and $a,b$  are the two invariants of the Maxwell field, 
related to $\bf E$, $\bf B$ by $a^2-b^2 = B^2-E^2,\quad ab = {\bf E}\cdot {\bf B}$.
This Lagrangian holds the information on the QED $N$ - photon amplitudes 
in the low-energy limit where all photon energies are small
compared to the electron mass, $\omega_i\ll m$. 
It corresponds to the diagrams shown in Fig. \ref{fig-EHL1loop}.

\vspace{1pt}

For the extraction of the amplitudes from the effective Lagrangian,
one expands it in powers of the Maxwell invariants,
\bear
{\cal L} (a,b) = \sum_{k,l} c_{kl}\, a^{2k}b^{2l} 
\ear
then fixes a helicity assignment and uses spinors helicity techniques \cite{56}.

\vspace{-5pt}


\section{Imaginary part of the effective action}

Except for the purely magnetic case where $b=0$, the proper-time integral in \eqref{ehspin} has poles on the 
integration contour at $ ebT =  k\pi$ which create an imaginary part.
For the purely electric case one gets \cite{schwinger51}
\begin{eqnarray}
{\rm Im} {\cal L}^{(1)}(E) &=&  \frac{m^4}{8\pi^3}
\beta^2\, \sum_{k=1}^\infty \frac{1}{k^2}
\,\exp\left[-\frac{\pi k}{\beta}\right]
\non\\
\label{schwinger}
\non
\end{eqnarray}
\vspace{-10pt}

($\beta = eE/m^2$). We note:

\begin{itemize}

\item
The $ k$th term relates to  coherent creation of $ k$  pairs in
one Compton volume. 


\item
$ {\rm Im}{\cal L}(E)$ depends on $ E$  non-perturbatively (non-analytically), which is consistent with
Sauter's \cite{sauter} interpretation of pair creation as  vacuum tunnelling (Fig. \ref{fig-pairtunnel}).

\vspace{1pt}

\begin{figure}[h]
\hspace{150pt}
{
\includegraphics[scale=0.5]{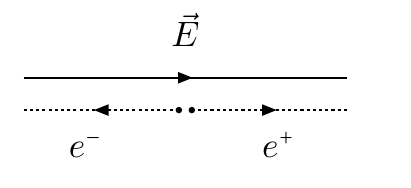}
}
\vspace{-5pt}

\caption{Pair creation by an external field as vacuum tunnelling.}
\label{fig-pairtunnel}
\end{figure}

\end{itemize}

\vspace{180pt}

\section{Beyond one loop}

\vspace{10pt}

\begin{figure}[h]
\hspace{120pt}
{
\includegraphics[scale=.7]{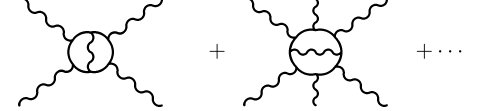}
}
\vspace{10pt}
\caption{Feynman diagrams contributing to the 2-loop EHL.}
\label{fig-2loopEHL}
\end{figure}

The two-loop correction to the EHL due to one internal photon exchange (Fig. \ref{fig-2loopEHL}) has been analyzed \cite{ritusspin,ditreu-book,18}, 
and turned out to contain important information on the Sauter tunnelling picture \cite{lebrit},
on-shell versus off-shell renormalization\cite{ritusspin,ritusscal}, and the asymptotic properties of the QED photon
S-matrix \cite{111}.

It leads to rather intractable two-parameter integrals. 
However, in the electric case its imaginary part $ {\rm Im}{\cal L}^{(2)}(E)$ permits a decomposition
analogous to Schwinger's \eqref{schwinger} \cite{lebrit}. For single-pair production, this is now interpreted as a tunnelling process where,
in the process of turning real, the electron-positron pair is already interacting at the one-photon exchange level.
%

Even for the imaginary part no completely explicit formulas are available. However, it
simplifies dramatically in the weak-field limit, where it just becomes an $\alpha\pi$ correction to the 
one-loop contribution:
 
 \bear
{\rm Im} {\cal L}^{(1)} (E) +
{\rm Im}{\cal L}^{(2)} (E) 
\,\,\,\, {\stackrel{\beta\to 0}{\sim}} \,\,\,\,
 \frac{m^4\beta^2}{8\pi^3}
\bigl(1+\alpha\pi\bigr)
\,{\rm e}^{-{\pi\over\beta}}
\non
\label{Im1plus2}
\ear
This suggests that higher loop orders might lead to an exponentiation, and indeed Lebedev and Ritus \cite{lebrit}
provided strong support for this hypothesis by showing that, assuming that
 \bear
{\rm Im} {\cal L}^{(1)} (E) +
{\rm Im}{\cal L}^{(2)} (E) 
+
{\rm Im}{\cal L}^{(3)} (E) + \ldots
\,\,\,\, {\stackrel{\beta\to 0}{\sim}} \,\,\,\,
 \frac{m^4\beta^2}{8\pi^3}
\,{\rm exp}\Bigl[ -{\pi\over\beta}+\alpha\pi \Bigr]
= {\rm Im}{\cal L}^{(1)}(E)\,\,{\rm e}^{\alpha\pi}
 \nonumber
\ear
then the result can be interpreted in the tunnelling picture as the  corrections to the
Schwinger pair creation rate due to the pair being created with a negative Coulomb
interaction energy

\begin{figure}[h]
\centering
\hspace{10pt}\includegraphics[scale=.5]{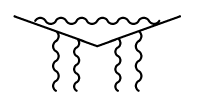}
\raisebox{1.1 em} {$\hspace{20pt}\Longleftrightarrow\hspace{20pt}$}
\includegraphics[scale=.4]{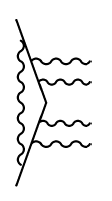}
\caption{Photon-corrected pair-creation vs. electron propagation in the field.}
\label{fig-crossed}
\end{figure}

\bear
m(E) \approx m + \delta m(E), \quad \delta m(E) =  - \frac{\alpha}{2}\frac{eE}{m} 
\nonumber
\ear
Moreover, the resulting field-dependent mass-shift $\delta m (E)$ is identical with 
the {\it Ritus mass shift}, originally derived by Ritus in \cite{ritusmass} from the crossed process of 
one-loop electron propagation in the field (Fig. \ref{fig-crossed}).

\vspace{5pt}

Unbeknownst to the authors of \cite{lebrit}, for scalar QED the corresponding conjecture had already been established 
two years earlier by Affleck, Alvarez and Manton \cite{afalma} using Feynman's  worldline path integral formalism 
and a semi-classical {\it worldline instanton} approximation.

%
%
%
%

\vspace{2pt}

\begin{figure}[h]
\hspace{70pt}
{
\includegraphics[scale=.6]{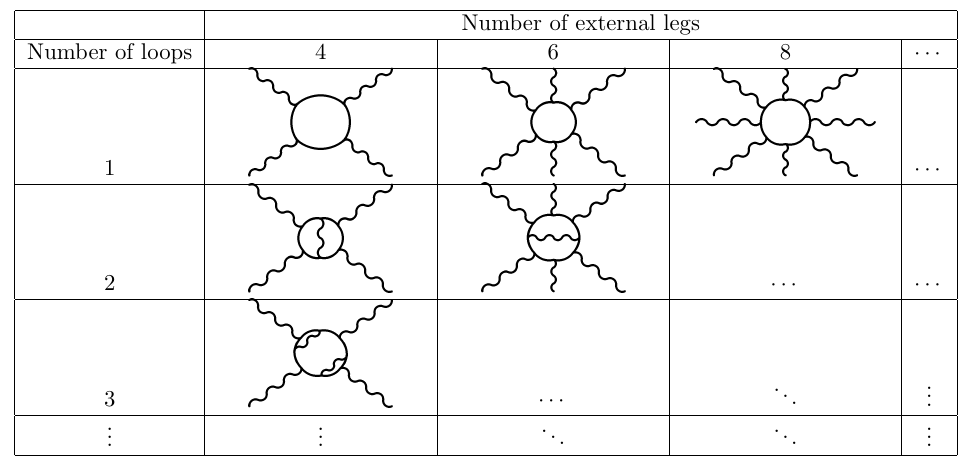}
}
\caption{Feynman diagrams contributing to the exponentiation hypothesis.}
\label{AAMfeyn}
\end{figure}

Diagrammatically, we note the following features of the exponentiation formula (see Fig. \ref{AAMfeyn}):

\begin{itemize}

\item
It Involves diagrams with any numbers of loops and legs.

\item
Although not shown, also all  the counter-diagrams from mass renormalization must contribute. 

\item
It does {\it not} include diagrams with more than one fermion loop 
(those get suppressed in the weak-field limit \cite{afalma}). 

\item
Horizontal summation produces the Schwinger exponential $ {\rm e}^{-\frac{\pi}{\beta}}$.

\item
Vertical summation produces the Ritus-Lebedev/Affleck-Alvarez-Manton exponential $ {\rm e}^{\alpha\pi}$.

\end{itemize}

\section{QED in 1+1 dimensions}


The exponentiation conjecture has so far been verified only at two loops. 
A three-loop check is in order, but calculating the three-loop EHL in $ D=4$ is presently hardly feasible. 
Motivated by work by Krasnansky \cite{Kras} on the EHL in various dimensions,
%
%
%
%
%
%
%
%
%
%
%
%
%
%
%
%
%
%
%
%
%
%
in 2010 two of the authors with D.G.C. McKeon started investigating the analogous problem in $2D$ QED.
In \cite{81} we used the worldline instanton method to generalize the exponentiation conjecture to the $2D$ case,
resulting in
\bear
{\rm Im}{\cal L}^{(all-loop)}_{2D}
\sim
\e^{-\frac{m^2\pi}{eE} + \tilde\alpha \pi^2  \kappa^2}
\label{exp2D}
\ear
where $\kappa =  m^2/(2ef)$, $ f^2=\fourth F_{\mu\nu}F^{\mu\nu}$,
and $\tilde\alpha = \frac{2e^2}{\pi m^2}$ is the two-dimensional analogue of the fine-structure constant. 
Defining the weak-field expansion coefficients in $2D$ by
\bear
{\cal L}^{(l)(2D)}(\kappa) &=& \frac{m^2}{2\pi} \sum_{n=1}^{\infty}(-1)^{l-1}c_{2D}^{(l)}(n) (i\kappa)^{-2n}
\label{defcn}
\ear
we then used Borel analysis to derive from (\ref{exp2D}) a formula for the limits of ratios of $l$ - loop to one - loop
coefficients:
\bear
{{\rm lim}_{n\to\infty}} {c^{(l)}_{2D}(n)\over c^{(1)}_{2D}(n+l-1)} 
&=& {(\tilde\alpha\pi^2)^{l-1}\over (l-1)!} 
\label{AAM2Dcoeff}
\ear
Moreover, we calculated the $2D$ EHL at one and two loops,
\bear
{\cal L}^{(1)}(f) &=& -{m^2\over 4\pi} {1\over\kappa}
\Bigl[{\rm ln}\Gamma(\kappa) - \kappa(\ln \kappa -1) +
\half \ln \bigl({\kappa\over 2\pi}\bigr)\Bigr]
\label{1loopEHL2D}\\
{\cal L}^{(2)}(f) &=& {m^2\over 4\pi}\frac{\tilde\alpha}{4}
\Bigl[ \tilde\psi(\kappa) + \kappa \tilde\psi'(\kappa)
+\ln(\lambda_0 m^2) + \gamma + 2 \Bigr]
\label{2loopEHL2D}
\ear \no
where $\tilde\psi (x) \equiv \psi(x) - \ln x + {1\over 2x}$, $\psi(x)=\Gamma^\prime(x)/\Gamma(x)$, and the
constant $\lambda_0$ comes from an IR cutoff. One finds from \eqref{1loopEHL2D} and \eqref{2loopEHL2D} that
\bear
c_{2D}^{(1)}(n) &=& (-1)^{n+1} \frac{B_{2n}}{4n(2n-1)} \\
c_{2D}^{(2)}(n) &=& (-1)^{n+1} \frac{\tilde\alpha}{8}\frac{2n-1}{2n}B_{2n}
\ear
Using properties of the Bernoulli numbers $ B_n$ it is then easy to verify that
\bear
\lim_{n\to\infty}  {c_{2D}^{(2)}(n)\over c_{2D}^{(1)}(n+1)} 
&=&
\tilde\alpha \pi^2
\nonumber
\ear
in accordance with \eqref{AAM2Dcoeff}.

\section{Three-loop EHL in 2D: diagrams}

\begin{figure}[h]
\hspace{140pt}
{
\includegraphics[scale=.4]{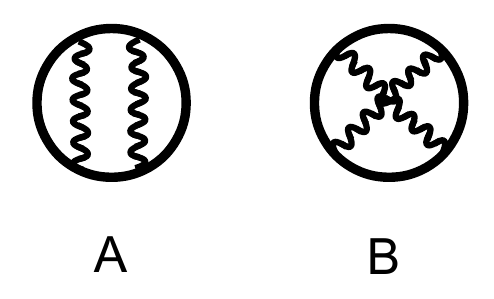}
}
\caption{Three-loop diagrams contributing to the exponentiation conjecture.}
\label{fig-AB}
\end{figure}

At three loops, we face the task of computing the two diagrams shown in Fig. \ref{fig-AB}
(there are also diagrams involving more than one fermion-loop, including several that
involve Gies-Karbstein tadpoles \cite{GK}, but those can be shown to be subdominant in the asymptotic limit).

The fermion propagators in these diagrams are the exact ones in the constant external field. Thus, although they are depicted as
vacuum diagrams, they are equivalent to the full set of ordinary diagrams of the given topology with any number of zero momentum
photons attached to them in all possible ways. 

Due to the super-renormalizability of $2D$ QED these diagrams are already UV finite.
They suffer from spurious IR - divergences, but those can be removed by going to the {\it traceless gauge} $\xi = -2$ \cite{121}. 
The calculation of diagram A is relatively straightforward, thus we focus on the much more substantial task of
computing diagram B and its weak-field expansion coefficients.

\begin{figure}[h]
\begin{center}
\includegraphics[scale=.7]{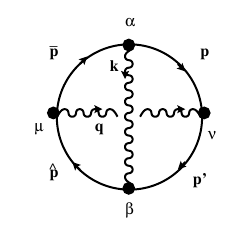}
\end{center}
\vspace{-10pt}
\caption{Parametrization of diagram B.}
\label{fig-DiagB}
\end{figure}

Introducing Schwinger parameters for this diagram as shown in Fig. \ref{fig-DiagB} leads to the integral representation \cite{121}
\begin{eqnarray}
{\cal L}^{3B} (f)
&=& \frac{\tilde\alpha^2m^2}{128\pi} \int_{0}^{\infty}dw dw' d\hat w d\bar w ~ I_B ~ \e^{-a} \nonumber\\
I_B &=& 
\frac{\rho^3}{\cosh^2 \rho w \cosh^2 \rho w' \cosh^2 \rho \hat w\cosh^2 \rho \bar w}
\frac{B}{A^3C} 
\nonumber\\&&
- \rho \frac{\cosh(\rho \tilde{w}) }{\cosh \rho w \cosh \rho w' \cosh \rho \hat w\cosh \rho \bar w}
\Bigl\lbrack
\frac{1}{A} - \frac{C}{G^2}\ln \Bigl(1+ \frac{G^2}{AC}\Bigr)\Bigr\rbrack
\nonumber
\label{Bw}
\end{eqnarray}
where
\bear
B &=& ({{\tanh}}^2 z + {\tanh}^2 \hat z)({\tanh} z' + {\tanh} \bar z) +  ({\tanh}^2 z' + {\tanh}^2 \bar z)({\tanh} z + {\tanh} \hat z) \non\\
C &=& {\tanh} z\,{\tanh} z' \,{\tanh}\, \hat z + {\tanh} z\, {\tanh} z'\, {\tanh} \bar z + {\tanh} z\, {\tanh} \hat z\, {\tanh} \bar z + {\tanh} z'\, {\tanh} \hat z \,{\tanh} \bar z \non\\
G &=& {\tanh} z\, {\tanh} \hat z - {\tanh} z'\, {\tanh} \bar z \nonumber
\ear
($ z = \rho w$  etc.).
Although for a three-loop diagram this is a fairly compact representation, an exact calculation is out of the question, and a straightforward expansion in powers of the external field
to get the weak-field expansion coefficients turns out to create huge numerator polynomials. To deal with those, we will now take advantage of the high symmetry
of the diagram. 

\section{Integration-by-parts algorithm}

Introduce the operator
$\tilde d  \equiv  \frac{\partial}{\partial w}- \frac{\partial}{\partial w'}+ \frac{\partial}{\partial \hat w} -\frac{\partial}{\partial \bar w} $
which acts simply on the trigonometric building blocks of the integrand. 
Integrating by parts with this operator, it is possible to write the integrand of $\beta_n$, the $n$-th coefficient of the expansion of $I_B$ as a power series in $\rho$, as a total derivative
$ \beta_n = \tilde d \theta_n$. Then, using once more the symmetry of the graph,
\bear
\int_{0}^{\infty}dw dw' d\hat w  d\bar w \,\e^{-a} \beta_n &=& \int_{0}^{\infty}dw d\bar w d\hat w dw'
 \tilde d
   \,\e^{-(w+w'+\hat w +\bar w)}
\theta_n 
\nonumber\\
&=&
4\int_{0}^{\infty}dw dw' d\hat w \,\e^{-(w+w'+\hat w)}\, \theta_n\vert_{\bar w =0}
\nonumber
\ear
The remaining threefold integrals are already of a fairly standard type. 

\section{Using the polynomial invariants of $D_4$}

Diagram $ B$ has the symmetries 
\bear
\qquad  w &\leftrightarrow& \hat w \nonumber\\
\quad w'  &\leftrightarrow& \bar w\nonumber\\
\quad (w,\hat w) &\leftrightarrow& (w',\bar w) \nonumber
\ear
Those generate the dihedral group $ D_4$. After a slight generalization to the inclusion of semi-invariants (invariants up to a sign) \cite{120},
Theorem 1 can be used to deduce that the numerator polynomials can be rewritten 
as polynomials in the variable $ \tilde w = w - w' +\hat w -\bar w$ with coefficients that are polynomials
in the four $ D_4$ - invariants $ a,v,j,h$,
\bear
a &=& w + w' + \hat w + \bar w \nonumber\\
v &=& 2 (w \hat w + w' \bar w) + (w + \hat w)(w' + \bar w) \nonumber\\
j &=& a \tilde w - 4 ( w \hat w - w' \bar w) \nonumber\\
h &=& a  (ww'\hat w + ww'\bar w + w\hat w \bar w + w'\hat w\bar w) + ( w \hat w - w' \bar w)^2 \nonumber
\ear
These invariants are moreover chosen such that they are  annihilated by $ \tilde d$.
Thus they are well-adapted to the integration-by-parts algorithm.
This rewriting leads to a very significant reduction in the size of the expressions  generated by the expansion in the field.

%
%
%
%
%
%
%
%
%
%
%
%
%
%
%
%
%
%
%
%
%
%
%
%
%
%
%
%

\section{Results}

In this way we obtained the first two coefficients of the weak-field expansion analytically,
\bear
\Gamma^B_0 &=& -\frac{3}{2}+\frac{7}{4}\zeta(3) \nonumber\\
\Gamma^B_1 &=& -\frac{251}{120} + \frac{35}{16}\zeta(3)\nonumber
\ear
and five more coefficients numerically (the coefficients $\Gamma_n$ are related to the ones introduced in
\eqref{defcn} by $4^n c^{(3)}_n = \frac{\tilde\alpha^2}{64}\Gamma_n$).
For a definite conclusion concerning the
exponentiation conjecture this is still insufficient, and the computation of further coefficients is in progress. 

\section{Outlook}

\begin{itemize}

\item
Writing Feynman graph polynomials in terms of invariant polynomials is a universal option that, to the best of our knowledge,
has not previously been used, but we expect that it will be found very useful for multiloop calculations involving diagrams 
with nontrivial symmetry groups and a large number of propagators. 

\item
In particular, this is the case for the weak-field expansion of the QED effective Lagrangian starting from three loops (in any dimension). 

\end{itemize}

\nolinenumbers

\end{document}